\documentclass{amsproc}
\newtheorem{theorem}{theorem}[section]
\newtheorem{lemma}[theorem]{Lemma}
\newtheorem{proposition}[theorem]{Proposition}
\newtheorem{corollary}[theorem]{Corollary}
\newtheorem{sublemma}[theorem]{Sublemma}

\theoremstyle{definition}
\newtheorem{definition}[theorem]{Definition}
\newtheorem{example}[theorem]{Example}

\newtheorem{notation}[theorem]{Notation}
\newtheorem{question}[theorem]{Question}

\theoremstyle{remark}
\newtheorem{remark}[theorem]{Remark}

\numberwithin{equation}{section}

\newcommand{\abs}[1]{\lvert#1\rvert}

\newcommand{\bd}{\noindent {\sc Proof}.\ \ }

\begin{document}

\title[The dichotomy of harmonic measures]
{The dichotomy of harmonic measures of compact hyperbolic
laminations}

\author{Shigenori Matsumoto}


\thanks{The author is partially supported by Grant-in-Aid for
Scientific Research (C) No.\ 20540096.}


\date{January 4, 2010}

\newcommand{\LL}{{\mathcal L}}
\newcommand{\R}{{\mathbb R}}
\newcommand{\DD}{{\mathbb D}^{d+1}}
\newcommand{\boundary}{{{\mathbb S}^d_\infty}}
\newcommand{\PP}{{\mathbb P}}
\newcommand{\GG}{{\mathbb G}}
\newcommand{\FF}{{\mathcal F}}
\newcommand{\EE}{{\mathbb E}}
\newcommand{\BB}{{\mathbb B}}
\newcommand{\HH}{{\mathcal H}}
\newcommand{\PPP}{{\mathcal P}}
\newcommand{\CC}{{\mathcal C}}
\newcommand{\oboundary}{{\mathbb S}^1_\infty}
\newcommand{\Q}{{\mathbb Q}}
\newcommand{\N}{{\mathbb N}}

\begin{abstract}
Given a harmonic measure $m$
of a hyperbolic lamination $\mathcal L$ on a compact metric
space $M$, 
a positive harmonic function $h$ on the universal cover
of a typical leaf is defined in such a way that the measure $m$ is
described in terms of these functions $h$ on various leaves. 
We discuss some properties of the function $h$. We show that
if $m$ is ergodic and not completely invariant, then
$h$ is typically unbounded and is induced by a probability $\mu$
of the sphere at infinity which is singular to the Lebesgue
measure. 
A harmonic measure is called Type I 
(resp.\ Type II) if
for any typical leaf, the measure $\mu$ is a point mass
(resp.\ of full support). We show that any ergodic
harmonic measure is either of type I or type II.
\end{abstract}

\maketitle

\section{Introduction}

We call $(M,\LL,g)$ a compact $C^2$ lamination if
$\LL$ is an $n$ dimensional lamination of class $C^2$ on a compact
metric space $M$ and if $g$ is a leafwise Riemannian metric of class
$C^2$.
(For the precise definition, see Sect.\ 2.)\ \
Then the leafwise Laplacian $\Delta f$ is defined for any continuous
leafwise $C^2$ function on $M$. A probability measure $m$ on $M$ is called
{\em harmonic} if for any such $f$, we have $m(\Delta f)
=0$.
A harmonic measure always exists for any compact $C^2$ lamination.

Given a harmonic measure $m$, there is a saturated conull set
$M^*$ such that a positive harmonic function $h$, called the
characteristic harmonic function, is defined
on the universal cover of each leaf in $M^*$ up to a constant
multiple. 
This function is obtained in the way of
describing the measure $m$ on each local chart. 
We first show (Theorem \ref{12}) that if $m$ is
ergodic and not completely invariant, then for any leaf in a saturated
conull set,
the characteristic harmonic function $h$ is unbounded.

A compact $C^2$ lamination $(M,\LL,g)$ of dimension $d$+1
is called {\em hyperbolic} if the metric $g$ has
curvature -1 on each leaf. The universal cover of each leaf is isometric
to the hyperbolic space $\DD$,
and the characteristic harmonic 
function $h$ corresponds to a probability
measure
$\mu$ on
the boundary at infinity $\boundary$. It depends upon
the choice of a base point in $\DD$, but its equivalence class is
uniquely determined by the leaf. 
We show (Theorem \ref{13}) that if $m$ is ergodic harmonic and not
completely invariant, then for any leaf in a saturated conull set,
the measure $\mu$ is
singular to the Lebesgue measure of $\boundary$.

\begin{definition}
A harmonic measure $m$ is called of {\em Type I} \ if the measure
 $\mu$ of $\boundary$ is a point mass for
any leaf in a saturated conull set, and of {\em Type II} \ if the support of 
$\mu$ is the whole $\boundary$.
\end{definition}

The main theorem of this paper is the following.

\begin{theorem} \label{1}
An ergodic harmonic measure is either of Type I or of Type II.
\end{theorem}

In Sect.\ 2, we prepare some necessary prerequisites about
harmonic measures. Especially the characteristic harmonic
function $h$ is defined. In Sect.\ 3, after a brief 
description of the
leafwise Brownian motion, we study its reverse process.
The reverse process plays a crucial role
in the proof of the unboundedness of $h$ (Sect.\ 3)
and the singularity of $\mu$ (Sect.\ 4).
In Sect.\ 5, we study the leafwise unit tangent
bundle $N$ of a compact hyperbolic lamination $\LL$.
There is a naturally defined lamination $\HH$ on $N$
of the same dimension as $\FF$.
Generalizing a result in [BM],
we discuss one to one correspondance between harmonic measures
on $\LL$ and pointed harmonic measures on $\HH$,
the latter being defined in Sect.\ 5.
Finally the proof of Theorem \ref{1}, as well as
some examples, is given in Sect.\ 6.

The author is grateful to Masahiko Kanai for supplying
him necessary knowledge about positive harmonic functions,
and to the referee for valuable comments which are helpful
for the improvement of the paper.

\section{Harmonic measure}

Let $M$ be a compact metric space, covered by a finite number
of open sets $E_i$. Assume there is a homeomorphism
$\varphi_i:E_i\to U_i\times Z_i$, where $U_i$ is an open ball in
$\R^n$ and $Z_i$ is a locally compact metric space.
 If $E_i\cap E_j\neq\emptyset$, then the transition function
$\psi_{ji}=\varphi_j\circ\varphi_i^{-1}$ is defined as a homeomorphism
from $\varphi_i(E_i\cap E_j)$ onto $\varphi_j(E_i\cap E_j)$.
Assume that the transition function is of the form
$$
\psi_{ji}(u,z)=(\alpha(u,z),\beta(z)),
$$
where $\alpha$ and $\beta$ are continuous, and $\alpha$ is
of class $C^2$ with respect to the first coordinate $u$ and its
first and second derivatives are continuous in $z$.
A subset of $M$ is called a {\em plaque}
if it is of the form $\varphi_i^{-1}(U_i\times z)$,
and a {\em transversal} if $\varphi_i^{-1}(u\times Z_i)$.
A maximal pathwise connected countable union of 
plaques is called a {\em leaf}. This gives birth to a decompositon
$\LL$ of $M$ into leaves, which is called a lamination of class $C^2$.

A leaf naturally has a structure of $n$ dimensional $C^2$ manifold.
A field of leafwise metric tensors is called a  
{\em leafwise Riemannian metric of class $C^2$} if
its leafwise derivatives up 
to order 2 (including order 0) are continuous on $M$.
In this paper the triplet $(M,\LL,g)$ is simply refered to
as a {\em compact $C^2$ lamination}.
By the compactness of $M$, each leaf of $\LL$ is complete
and of bounded geometry.
The leafwise volume defined by $g$ is denoted by ${\rm vol}$.

Henceforce we depress the homeomorphism $\varphi_i$ and
consider $U_i\times Z_i$ as an open subset of $M$,
which is called a {\em local chart}.

A function $f:M\to\R$ is said to be {\em continuous leafwise $C^2$}
if it is of class $C^2$ in each leaf and its derivative up to
order 2 is continuous on $M$. Then the leafwise Laplacian
$\Delta f$ with respect to $g$ is defined, and is a continuous 
function on $M$.

\begin{definition}
A probability measure $m$ on $M$ is called {\em harmonic}
if $m(\Delta f)=0$ for any continuous leafwise $C^2$
function.
\end{definition}

\begin{remark}
A harmonic measure always exits for any compact
$C^2$ lamination $(M,\LL,g)$.
\end{remark}

See [C] Theorem 3.5 for a simple proof using the Hahn-Banach theorem.

Here is a structure theorem of a harmonic measure on a local chart.

\begin{theorem} \label{2}
Assume $m$ is a harmonic measure on a compact $C^2$ lamination.
For any local chart $U\times Z$, there are a measure $\nu$
on $Z$ and a
function $h:U\times Z\to\R$ 
with the following properties. 

{\rm (1)} $h$ is positive and $m$-measurable.

{\rm (2)} For $\nu$-a.a.\ $z$, the restriction of 
$h$ to the plaque $U\times z$
is harmonic and $h$vol is a probablity measure of the plaque.

{\rm (3)} For any continuous function with 
support in $U\times Z$, we have
$$
m(f)=\int_Z\int_{U\times{z}}f(u,z)h(u,z)d{\rm vol}(u)d\nu(z).
$$

Furthermore if a probability measure $m$ on $M$ is represented
in this way in any local chart, then $m$ is  harmonic.

\end{theorem}

\begin{notation}
The theorem says that the measure $m$ restricted to $U\times Z$
is disintegrated in such a way that the conditional probability measure on each fiber
$U\times z$ is $h(\cdot,z){\rm vol}$ and the push forward measure on the base
$Z$ is $\nu$. Henceforth this is denoted as
\begin{equation} \label{e}
m\vert_{U\times Z}=\int_Zh{\rm vol}\,d\nu.
\end{equation}
\end{notation}
 
\bd
By the disintegration theorem 
we have
$$m\vert_{U\times Z}=\int_{Z} m_z\,d\nu(z),
$$
where $m_z$ is a probability measure
on $U\times z$  and the assignment $z\mapsto m_z$ is measurable.
The measure $\nu$ is 
the push forward of $m$ by the projection $p_2:U\times Z\to Z$,
and is not necessarily a probability measure.

Denote the other projection by $p_1:U\times Z\to U$. 
The leafwise Riemannian metric on each plaque $U\times z$ is
transfered to a Riemannian metric on $U$, 
and the corresponding Laplacian on $U$ is denoted by $\Delta_z$.
Consider any function $f$ from the space $C^2_c(U)$ of
the $C^2$ functions with compact support,
and any continuous 
function $g$ on $Z$ with compact support.
Then the product $f\circ p_1\,g\circ p_2$
is a continuous leafwise $C^2$ function
whose support is contained in $U\times Z$ and we have 
$$\Delta(g\circ p_2\,f\circ p_1)=g\circ p_2\,\, \Delta(f\circ p_1)
\ \mbox{ and }\
m_z(\Delta(f\circ p_1))=m_z(\Delta_zf).
$$
Since $m(\Delta(g\circ p_2\,f\circ p_1))=0$, we have
$$\int_{Z}m_z(\Delta_z f)g(z)d\nu(z)=0.
$$
By the measurablility of the assignment $z\mapsto m_z$
and the boundedness of $\Delta_zf$,  $m_z(\Delta_z f)$ is an integrable function
on $Z$ and thus $m_z(\Delta_z f)\nu$ is a signed measure on $Z$,
for which an arbitrary compactly supported continuous
function $g$ integrates to 0.
This implies that for $\nu$-a.a.\ $z$, $m_z(\Delta_z f)=0$.

Since $C^2_c(U)$ has a countable dense subset $S$, there is
a $\nu$-conull set $Z^*$ of $Z$ such that if $z\in Z^*$,
 $m_z(\Delta_z f)=0$ for
any $f\in S$, and therefore for any $f\in C^2_c(U)$.
But as is well known ([N]), $m_z(\Delta_z f)=0$ for any 
$f\in C^2_c(U)$ if and only if $m_z=h_z\,{\rm vol}$ for a harmonic
function $h_z$ on $U$ with respect to the Laplacian $\Delta_z$. 
Setting $h(u,z)=h_z(u)$, we obtain (\ref{e}).
 
Next we are going to show that the function $h$
 is measurable.
Consider another measure $m_0$ on $U\times Z$, given
by $\int_Z {\rm vol}/{\rm vol}(U\times z)d\nu(z)$. Clearly $m$ and
$m_0$ are mutually equivalent and thus we have
$$
m=km_0$$ for some $m_0$-measurable (equivalently $m$-measurable) function $k$.
But the uniqueness of the disintegration implies that for $\nu$-a.a.\ $z$,
$$
h(u,z)=k(u,z)/{\rm vol}(U\times z),
$$
showing that $h$ is measurable.

Finally the converse statement is easy to show.
\qed

\medskip
As an immediate corollary, we have:

\begin{corollary} \label{2A}
If a function $f$ on $M$ is $C^2$ on each leaf and $\Delta f$
is $m$-integrable, then $m(\Delta f)=0$.
\end{corollary}

\begin{remark} In [G], harmonic measures are defined
by the property in Corollary \ref{2A}, and the structure
theorem is obtained. Our proof of Theorem \ref{2} shows the equivalence of
the two definitions.
\end{remark}

\smallskip
Suppose two local charts $U\times Z$ and $U'\times Z'$ intersect
and the harmonic measure $m$ is decomposed on each local chart
as 
$$m\vert_{U\times Z}=\int_Zh\,{\rm vol}\,d\nu\ \mbox{ and }\
m\vert_{U'\times Z'}=\int_{Z'}h'\,{\rm vol}\,d\nu',
$$
 then 
in the intersection of $\nu$-a.a.\ plaque $U\times z$
and $\nu'$-a.a.\ plaque $U'\times z'$, we have 
\begin{equation} \label{e1}
h'/h=d\nu/d\beta\nu',
\end{equation}
where $\beta$ is the holonomy map from (a part of) $Z'$ to $Z$.
On one hand this shows that $\nu$ and $\nu'$ are equivalent via
the holonomy map. More generally we have:

\begin{proposition} \label{3}
A harmonic measure $m$ is leafwise smooth, i.e.,

{\rm (1)} If a Borel set $B\subset M$ satisfies ${\rm vol}(B\cap L)=0$ for any leaf
$L$, then $m(B)=0$.

{\rm (2)} If a Borel set $B$ is $m$-null, then the set $\widehat B$ is also
$m$-null, where $\widehat B$ is the union of the leaves $L$ such
that ${\rm vol}(B\cap L)>0$.
\end{proposition}

On the other hand the equality (\ref{e1}) 
shows that on the intersection of two plaques,
the function $h'$ is a positive constant multiple of $h$.
Dividing $h'$ by that constant,
one can continuate $h$ along a chain of plaques. Of course
this does not yield a function on a leaf, since
there will be a monodromy for $h$. However we will get
a function on the holonomy cover of a leaf.

In what follows, when we say ``for $m$-a.a.\ leaf $L$'', 
this means `` there exists a saturated conull set 
$M^*$ and for any leaf $L$ in $M^*$''.

\begin{proposition} \label{4}
{\rm (1)} For $m$-a.a.\ leaf $L$, the function $h$ has
a well defined prolongation as a positive harmonic function
on the holonomy cover $\hat L$. On $\hat L$ two such functions
which start from different plaques are unique up to
a positive constant multiple.

{\rm (2)} Given a path in $L$ the ratio
of $h$ at the initial point and the terminal point of 
any lift of the path to $\hat L$ is the same.
\end{proposition}

\noindent
{\sc Proof}.
To see (1), cover $M$ with a finite number
of local charts $U_i\times Z_i$.  Then there is a $\nu$-conull set
$Z_i^*$ of each $Z_i$ such that the harmonic function $h$ is
defined on $U_i\times Z_i^*$. The saturation of the union of 
$Z_i\setminus Z_i^*$ is $m$-null by Proposition \ref{3},
and for any leaf $L$ in the complement $M^*$,
 the function $h$ has a  prolongation on its holonomy cover
$\hat L$.

The uniqueness part in (1), as well as
(2), follows immediately from the  construction.
\qed

\bigskip

Of course the harmonic function $h$ has a lift to the universal
covering space $\tilde L$ of  $m$-a.a.\ leaf $L$, which will be
denoted by the same letter $h$.
The above statement (2) holds also for lifts of paths to the universal
covering
space. Let $\Gamma$ be the deck transformation group of the covering
map $\tilde L\to L$. Then we have:

\begin{corollary} \label{4+}
For any $\gamma\in\Gamma$, $h\circ\gamma$ is a constant multiple
of $h$.
\end{corollary}

\noindent
{\sc Proof}. Join two points $x\in\tilde L$ and $y\in\tilde L$ 
by an arc $c$.
Then $\gamma x$ and $\gamma y$ are joined by $\gamma c$.
Two arcs $c$ and $\gamma c$ are lifts of the same arc in $L$.
Therefore we have
$$
h(y)/h(x)=h(\gamma x)/h(\gamma y).$$
Since $x$ and $y$ are arbitrary, this shows that the function
$h(\gamma x)/h(x)$ is independent of $x$.
\qed

\medskip

\begin{definition} The function $h$ in Proposition \ref{4} is called
the {\em characteristic harmonic function of $m$}.
\end{definition}

Notice that the characteristic harmonic function is defined only
up to a positive constant multiple.

\medskip
A harmonic measure $m$ is called {\em completely invariant} if
the characteristic harmonic functions are constant on
(the holonomy cover of) $m$-a.a.\ leaves. In this case
$m$ corresponds to a transverse invariant measure, i.e., 
an assignment of a finite measure to each transversal
which is invariant by the holonomy maps.
Conversely a transverse invariant measure gives rise to
a harmonic measure $m$ whose characteristic harmonic function is constant on
$m$-a.a.\ leaf.
Only a special class of laminations admit completely invariant
measures.

\section{Brownian motion and its reverse process}

Let $(M,\LL, g)$ be a compact $C^2$ lamination. Denote by
$\Omega$ the space of all the continuous leafwise paths 
$\omega:[0,\infty)\to M$ and for $t\geq0$, define a map
$X_t:\Omega\to M$ by $X_t(\omega)=\omega(t)$.
Let $\mathcal B$ be the $\sigma$-algebra of the Borel subsets
of $\Omega$ with respect to the compact open topology. It is well known,
easy to show, that $\mathcal B$ coincides with the minimal
$\sigma$-algebra for which $X_t$ ($t\geq 0$) is Borel.
In other words $\mathcal B$ is generated by the {\em cylinder sets}
$\{X_{t_1}\in B_1,\ldots,X_{t_r}\in B_r\}$ ($0\leq t_1<\cdots< t_r$,
$B_i$; a Borel subset of $M$).

The leafwise Riemannian metric $g$ gives the heat kernel
$p_t(x,y)$ ($t>0$) on each leaf. Define $p_t(x,y)$ 
for any two points $x,y\in M$ by setting $p_t(x,y)=0$
unless $x$ and $y$ lie on the same leaf. The heat kernel defines the
Wiener probability measure $W^x$ on $\Omega$ ($x\in M$).
For a cylinder set $\{X_{t_1}\in B_1,\ldots,X_{t_r}\in B_r\}$
($t_1>0$), we define
\begin{eqnarray}\label{eqn}
& W^x\{X_{t_1}\in B_1,\ldots, X_{t_r}\in B_r\}
\\
&=\int_{B_1}\cdots\int_{B_r}p_{t_1}(x,y_1)p_{t_2-t_1}(y_1,y_2)
\cdots p_{t_{r}-t_{r-1}}(y_{r-1},y_r)\,d{\rm vol}(y_r)\cdots d{\rm vol}(y_1).
\nonumber
\end{eqnarray}
Then $W^x$ satisfies the dropping condition, and
therefore it is defined not only for a cylinder set but also for any
set in $\mathcal B$. That is, $W^x$ is a probability measure on
$(\Omega,\mathcal B)$. It is concentrated on the subset
$\Omega_x=X_0^{-1}(x)$, since the probability measure 
$p_t(x,\cdot){\rm vol}$ tends to the Dirac mass at $x$ as $t\to0$. 

\begin{lemma} \label{24}
The system of measures $\{W^x\}_{x\in M}$ is Borel in the sense
that for any $S\in\mathcal B$, the assignment $x\mapsto W^x(S)$
is Borel. 
\end{lemma}

\bd
Let $\mathcal C$ be the family of the subsets $S$ in $\Omega$
for which $M\ni x \mapsto W^x(S)$ is Borel,
and let $\mathcal A_0$ be the finite algebra formed by 
finite disjoint unions of cylinder sets.
Then 
$\mathcal A_0$ is contained in $\mathcal C$.
For $\{X_{t_1}\in B_{t_1}\}\in\mathcal A_0$, see [CC] Lemma 2.3.1.
General case follows easily from this. 

For an isolated ordinal $\alpha>0$, define 
$\mathcal A_\alpha$ to be the family of a subset which is
obtained from subsets of $\mathcal
A_{\alpha-1}$,
by a finite sequence of two operations; one, 
taking a countable increasing union
and the other, countable decreasing intersection.
Then it is easy to show that $\mathcal A_\alpha$ forms a finite algebra.
Moreover $\mathcal A_\alpha$ is contained in $\mathcal C$,
since a pointwise limit of Borel functions is Borel.
For a limit ordinal $\alpha$, let 
$\mathcal A_\alpha=\bigcup_{\beta<\alpha}A_\beta$.
Then again $\mathcal A_\alpha$ is a finite algebra contained in $\mathcal C$.

The increasing sequence $\{\mathcal A_\alpha\}$ stabilizes.
Define $\mathcal A=\mathcal A_{\alpha_0}$, where
$A_\beta=A_{\alpha_0}$ for any
ordinal $\beta\geq \alpha_0$. Then $\mathcal A$ is contained in $\mathcal C$.
On the other hand $\mathcal A$
is  clearly a $\sigma$-algebra. 
Therefore any Borel set, an element of the minimal $\sigma$-algebra 
which contains $\mathcal A_0$, belongs to $\mathcal A$,
and hence to $\mathcal C$.
\qed

\bigskip

The expectation
of $W^x$ is denoted by $E^x$. 
Applying Lemma \ref{24}, one can show that for any bounded Borel 
function $f:M\to\R$, its diffusion
$D_tf$ is bounded Borel, where
$$
(D_tf)(x)=E^x[f(X_t)]=\int_M p_t(x,y)f(y)d{\rm vol}(y).
$$
More generally the diffusion operator $D_t$ defines a semigroup of contractions
on the space $L^p(M,m)$ ($1\leq p<\infty$) for
a harmonic measure $m$ and on $C(M)$,
the space of continuous functions ([C]). 

Since $\{W^x\}$ is a Borel system of measures,
by integrating $W^x$ over any probability measure $m$ on $M$,
we get a probability measure $\PP_m$ on $\Omega$, i.e., 
$$\PP_m=\int_MW^xdm(x).
$$ Precisely
for any bounded Borel function $F:\Omega\to\R$, 
$$
\PP_m(F)=\int_ME^x[F]dm(x).
$$

For $t\geq0$ let $\theta_t:\Omega\to\Omega$ denote the shift map by $t$, 
i.e., 
$(\theta_t\omega)(t')=\omega(t+t')$. 

\begin{theorem} \label{4A}
The probability measure $m$ is harmonic if and only if the probability measure
$\PP_m$ is $\theta_t$-invariant for any $t\geq0$.
\end{theorem}

For the proof, see [CC] Theorem 2.3.7.

A harmonic measure $m$ is called {\em ergodic} if whenever it is written as 
a nontrivial linear combination of two harmonic measures $m_1$
and $m_2$, we have $m=m_1=m_2$.

\begin{theorem} \label{5}
Let $m$ be a harmonic measure. Then the following conditions
are equivalent.

{\rm (1)} $m$ is ergodic.

{\rm (2)} For any saturated Borel set $M'$ in $M$, we have either $m(M')=0$ or
$m(M')=1$.

{\rm (3)} If $f\in L^1(M,m)$ satisfies $D_tf=f$ for any $t\geq0$, then $f$ is a
constant.

{\rm (4)} $\PP_m$ is ergodic with respect to the semiflow
defined by the shift map $\theta_t$, i.e., 
if a Borel subset $S$ satisfies $\theta_t^{-1}(S)=S$
for any $t\geq0$, then either $\PP_m(S)=0$ or $\PP_m(S)=1$.
\end{theorem}

That (1) $\Rightarrow$ (2) follows from Corollary \ref{2A},
and that (4) $\Rightarrow$ (1) is immediate. The other 
implications (2) $\Rightarrow$ (3) $\Rightarrow$ (4)
can be shown in exacly the same way as the proof of
Theorem 3.1 of [F].

\bigskip

The diffusion operator $D_t:L^2(M,m)\to L^2(M,m)$ ($m$ a harmonic measure)
is not self-adjoint unless $m$ is completely invariant.
Its adjoint $D_t^*$ is first considered in [K].
Let $h$ be the characteristic harmonic function defined on the 
holonomy cover $\hat L$ of $m$-a.a.\ leaf $L$.
Denote by $\hat p_t$ the heat kernel on $\hat L$. We have
$$
p_t(x,y)=\sum_{\hat y}\hat p_t(\hat x, \hat y),
$$
where the sum is taken for all the points $\hat y$
over $y$, and independent of the choice of $\hat x$ over $x$.

We shall summerize well known properties of the heat kernel $\hat p_t$
on $\hat L$ which follows from the bounded geometry of $\hat L$.

\begin{lemma}\label{last}
 For any harmonic function $g:\hat L\to\R$, we have
\begin{eqnarray*}
&g(\hat x)=\int_{\hat L}g(\hat y)\hat p_t(\hat x,\hat y)d{\rm vol}(\hat y)
\ \mbox{ and }\\
&\hat p_{t+t'}(\hat x,\hat z)=\int_{\hat L}
\hat p_t(\hat x,\hat y)\hat p_{t'}(\hat y,\hat z)d{\rm vol}(\hat y).
\end{eqnarray*}
\end{lemma}

Now define a new heat kernel on $\hat L$ by
$$
\hat q_t(\hat x,\hat y)
=\frac{h(\hat y)}{h(\hat x)}\hat p_t(
\hat x,\hat y).
$$

The following lemma follows immediately from Lemma \ref{last}.

\begin{lemma} \label{6}
We have
\begin{eqnarray} \label{e3}
&\int_{\hat L}\hat q_t(\hat x,\hat y)\,d{\rm vol}(\hat y)=1
\ \mbox{ and }\\
&\hat q_{t+t'}(\hat x,\hat z)=\int_{\hat L}
\hat q_t(\hat x,\hat y)\hat q_{t'}(\hat y,\hat z) \label{e4}
\,d{\rm vol}(\hat y).
\end{eqnarray}
\end{lemma}

Define a heat kernel $q_t$ on the leaf $L$ by
$$
q_t(x,y)=\sum_{\hat y}\hat q_t(\hat x, \hat y).
$$

\begin{theorem} \label{7}
The dual operator $D_t^*$ is expressed for any $f\in L^2(M,m)$ as
$$
(D_t^*f)(x)=\int_Lq_t(x,y)f(y)\,d{\rm vol}(y),
$$
where $L$ is the leaf through $x$.
\end{theorem}

Although this theorem is known to Vadim Kaimanovich, we shall include
a proof, since there seems to be none in the literature.

Let $\GG$ denote the holonomy groupoid associated to the lamination
$\LL$, i.e.,  $\GG$ is the space of leafwise paths modulo same end points
and identical holonomy germs. Denote by $r,\,s:\GG\to M$ the range 
and the source maps. The fiber $s^{-1}(x)$ is homeomorphic to
the holonomy cover of the leaf through $x$, and the corresponging volume form
of $s^{-1}(x)$ is denoted by ${\rm vol}_x$. Integrating these forms
(seen as measures) over the harmonic measure $m$ of $M$,
we get a measure $m_\GG$ on $\GG$. That is,
$$m_\GG=\int_M {\rm vol}_x dm(x).$$
 Likewise we define a measure ${\rm vol}^y$ on
$r^{-1}(y)$. Define a function $\varphi:\GG\to\R$ by 
$\varphi([\gamma])=h(\gamma(1))/h(\gamma(0))$, where $h$ is the
characteristic harmonic function which is defined on the holonomy
cover of $m$-a.a.\ leaf. The function $\varphi$ is well defined by
Proposition \ref{4} (2). Denote by $J:\GG\to\GG$ the inverse map.

\begin{lemma} 
We have $Jm_\GG=\varphi\cdot m_\GG$.
\end{lemma}

\bd
For an arbitrary $[\gamma]\in\GG$,
choose a neighbourhood $U\times V\times Z$ of $[\gamma]$ in $\GG$, where
$U\times Z$ is a local chart containing $\gamma(0)$ 
so that the holonomy along $\gamma$ is defined on $Z$ and $V$
is a leafwise neighbourhood of $\gamma(1)$. 
Changing the notations slightly, we consider $U$ (resp.\ $V$)
to be a neighbourhood of $\tilde \gamma(0)$ (resp.\ $\tilde \gamma(1)$)
in the universal cover $\tilde L$ of the leaf $L$, where $\tilde \gamma$
is a lift of $\gamma$
to $\tilde L$. 
Choosing $Z$ smaller if necessary, we may assume
that there is a precompact simply 
connected open
set $W$ of $\tilde L$ such that $U\cup V\cup\tilde\gamma\subset W$ and that
there is a lamination preserving embedding of $W\times Z$ into $M$.

Then by Theorem \ref{2},
$$m\vert_{W\times Z}=\int_Zh\,{\rm vol}\,d\nu
$$ 
for a leafwise
harmonic function $h$ and a measure $\nu$ on $Z$.
For $(u,v,z)\in U\times V\times Z\subset\GG$,
denote $s(u,v,z)=(u,z)=x$ and $r(u,v,z)=(v,z)=y$.
Restricted to $U\times V\times Z$, ${\rm vol}_x$ is
the volume form on $u\times V\times z$ and ${\rm vol}^y$
on $U\times v\times z$.

On $U\times V\times Z$ we have
$$
m_\GG=\int_Z {\rm vol}_x\cdot h(x){\rm vol}^y\,d\nu.$$
On the other hand on $V\times U\times Z$, 
$$
Jm_\GG=\int_Z  {\rm vol}_y\cdot h(y){\rm vol}^x\,d\nu
=\int_Z \varphi\cdot {\rm vol}_y\cdot h(x){\rm vol}^x\,d\nu
=\varphi\cdot m_{\GG},$$
showing the lemma.
\qed

\smallskip

\begin{remark}
The measure $m_\GG$ is defined not only for a harmonic measure, but
also for any probability measure $m$ on $M$. It is interesting
to remark that the leafwise smoothness (Proposition \ref{3}) of $m$ is 
equivalent to a basic notion in measured groupoids, 
the equivalence of $Jm_\GG$ with $m_\GG$ [AR].
\end{remark}

\smallskip
\noindent
{\sc Proof of  Theorem \ref{7}.}
The Riemannian heat kernel on the holonomy cover of the leaf yields
a function $\check p_t$ on $\GG$ by
$$\check p_t([\gamma])=\hat p_t(\gamma(0),\gamma(1)).$$
Notice that $\check p_t\circ J=\check p_t$.
 Likewise
a function $\check q_t$ is defined from $\hat q_t$. They satisfy
$\check q_t=\varphi\check p_t$.
Clearly we have
$$
(D_tf)(x)=\int_{s^{-1}(x)}\check p_t\, f\circ r \, d{\rm vol}_x.$$
Thus
\begin{eqnarray*}
&\langle  D_tf,g\rangle=\int_M(\int_{s^{-1}(x)}
\check p_t\, f\circ r \, d{\rm vol}_x)g(x)dm(x)
=
\int_\GG \check p_t\,f\circ r\,g\circ s\, dm_\GG\\
&=\int_\GG \check p_t\, f\circ s\, g\circ r\, \varphi\,dm_\GG
=\int_\GG \check q_t\, g\circ r\, f\circ s\, dm_\GG\\
&=\int_M(\int_{s^{-1}(x)}\check q_t\, g\circ r\, d{\rm vol}_x)
f(x)dm(x)
=\langle f,D_t^*g\rangle.
\end{eqnarray*}
Therefore we have
$$
(D_t^*g)(x)=\int_{s^{-1}(x)}\check q_t\, g\circ r\, d{\rm vol}_x
=\int_L q_t(x,y)g(y)\, d{\rm vol}(y),$$
completing the proof.
\qed

\medskip

Now let us define the reverse process. 
First of all extend the new heat kernel $q_t$ to $M\times M$,
by putting $q_t(x,y)=0$ unless $x$ and $y$ lie on the
same leaf.
Let $\Omega_-$ be the
space of continuous leafwise paths $\omega$ from $(-\infty,0]$ to $M$,
with the random variable $X_{-t}:\Omega_-\to M$ defined
by $X_{-t}(\omega)=\omega(-t)$ ($t\geq0$).
For $x\in M$, define the Wiener measure $W^x_-$ on $\Omega_-$
using the kernel $q_t$, that is, for example for 
$0<t_1< t_2$ and for any Borel sets $B_1$ and $B_2$ of $M$,
$$
W^x_-\{X_{-t_2}\in B_2, X_{-t_1}\in B_1\}
=\int_{B_2}\int_{B_1} q_{t_1}(x,y)q_{t_2-t_1}(y,z)\,d{\rm vol}(y)\,d{\rm vol}(z).$$

Lemma \ref{6} implies that $W^x_-$ is a probability 
measure, a probability because of
(\ref{e3}), the dropping condition guaranteed by (\ref{e4}).
The kernel $q_t$ clearly satisfy the normal estimate of Cheng, Li
and Yau ([CLY]) since the ratio to the Riemannian heat kernel is
controlled by the Harnack inequality;
the logarithm of any positive harmonic
function defined on the holonomy cover of any leaf of $\LL$ is uniformly 
Lipschitz (due to the uniform boundedness of geometry of leaves). 
Therefore the reverse Wiener measure
$W^x_-$ is concentrated on the set of continuous paths.
Moreover it is concentrated on the subspace $\Omega_{-,x}=X_{0}^{-1}(x)$.

Now let $\overline \Omega$ be the space of biinfinite continuous leafwise paths
$\omega:\R\to M$. Denote the like defined random variable 
by the same letter $X_t:\overline\Omega\to M$ for $t\in\R$.
Also denote $\overline\Omega_x=X_0^{-1}(x)$.
Then by the natural identification of $\Omega_{-,x}\times\Omega_{x}$ 
with $\overline\Omega_x$,
the product measure $W_-^x\times W^x$ is considered to be
a measure on $\overline \Omega_x$, or on $\overline \Omega$.

Define a probability measure
$\overline\PP_m$ on 
$\overline \Omega$
by 
$$\overline \PP_m=\int_MW^x_-\times W^x\, dm(x).
$$
Denote its
expectation by $\overline\EE_m$.
Let $\theta_t:\overline\Omega\to\overline\Omega$ be the shift map.

\begin{proposition} \label{8}
The shift map $\theta_t:\overline\Omega\to\overline\Omega$ preserves the measure
$\overline\PP_m$.
\end{proposition}

\bd
We shall raise one example of computation.
$$
\overline\PP_m\{X_{-t}\in B, X_{t'}\in B'\}
=\int_M dm(x)\int_B q_t(x,y)dy\int_{B'}p_{t'}(x,z)dz
$$
$$
=\langle D_t^*\chi_B, D_{t'}\chi_{B'}\rangle_m
=\langle \chi_B, D_tD_{t'}\chi_{B'}\rangle
=\overline \PP_m\{X_0\in B, X_{t+t'}\in B'\}.$$
\qed

\begin{theorem} \label{9}
If $m$ is an ergodic harmonic measure, then $\overline \PP_m$ is
ergodic with respect to the flow $\{\theta_t\}$.
\end{theorem}

Before starting the proof, we
recall the definition of conditional expectations.
Denote by $\overline\FF$ the $\sigma$-algebra formed by
the $\overline\PP_m$-measurable subsets.
For $t\in\R$, let $\overline\FF_t$ be the minimal complete $\sigma$-algebra
for which the map $X_s$ is measurable for any $s\geq t$.

For example, in order to understand $\overline\FF_0$,
consider the measurable partition of $\overline\Omega$
defined by the natural projection 
$\pi:\overline\Omega\to\Omega$.
Then $\overline\FF_0$ consists of measurable subsets saturated
by this partition.
A function $F$ is $\overline\FF_0$-measurable if
and only if there is a measurable function $H$ on $\Omega$
such that $F=H\circ\pi$.

For any integrable function $F:\overline\Omega\to\R$,
denote by $\overline\EE_m\,[F\mid\overline\FF_t]$ the conditional expectation
with respect to $\overline\FF_t$. This is a unique $\overline\FF_t$-measurable 
function on $\overline\Omega$ such that for any
bounded $\overline\FF_t$-measurable
function $G$,
$$\overline\EE_m\,[G\,\overline\EE_m\,[F\mid\overline\FF_t]\,]=
\overline\EE_m\,[GF].
$$
 
One word of explanation for the geometer readers.
$\overline\FF_t$ defines a measurable partition of $\overline\Omega$:
almost all classes of the partition admit the conditional probability measure. 
Integrating $F$ by the conditional probability measure we obtain a measurable function
on the quotient space. But it is customary, more convenient, to
consider it to be a
$\overline\FF_t$-measurable function $\overline\EE_m\,[F\mid\overline\FF_t]$
defined on the total space $\overline \Omega$.

\smallskip
\noindent
{\sc Proof of Theorem \ref{9}.}
For an integrable function $F$ on $\overline\Omega$ define the
Birkhoff average $\BB F$ by
$$
\BB F=\lim_{t\to\infty}\frac{1}{t}\int_0^tF\circ \theta_s ds.
$$
By the ergodic theorem, the operator $\BB$ is a well defined contraction
on $L^1(\overline\Omega,\overline\PP_m)$, which is $\theta_t$-invariant.
 
Since by Theorem \ref{5}, $\theta_t$ is ergodic in $(\Omega,\PP_m)$,
the Birkhoff average $\BB F$ is constant if $F$ is $\overline\FF_0$-measurable.
Moreover this holds for any $\overline\FF_{-t}$-measurable function $F$
for any $t$, since then $F\circ \theta_{t}$ is $\overline\FF_0$-measurable
and $\BB F=\BB(F\circ\theta_{t})$.

For any bounded $\overline \FF$-measurable function $F$,
the $\overline\FF_{-n}$-measurable function 
$F_{-n}=\overline\EE_m\,[F\mid \overline\FF_{-n}]$ converges to $F$ pointwise,
by the martingale convergence theorem ([O] Appendix C).
Thus we have $\BB F_{-n}\to \BB F$, and since $\BB F_{-n}$ is constant,
the function $\BB F$ is also constant, showing the ergodicity.
\qed

\bigskip

Applying the Birkhoff theorem to $f\circ X_0:\overline \Omega\to\R$ 
for a continuous function $f:M\to \R$ by virtue of Theorem \ref{9},
we have $\overline\PP_m$-almost surely
$$
\lim_{t\to\infty}\frac{1}{t}\int_0^tf(X_s) ds
=\lim_{t\to\infty}\frac{1}{t}\int_{-t}^0f(X_s) ds 
=m(f).
$$

Equivalently, denoting the Dirac mass by $\delta_.$, we have 
$\overline\PP_m$-almost surely
\begin{equation} \label{e21}
\lim_{t\to\infty}\frac{1}{t}\int_0^t\delta_{X_s} ds
=\lim_{t\to\infty}\frac{1}{t}\int_{-t}^0\delta_{X_s} ds 
=m,
\end{equation}
where the limit is taken in the space of the probability
measures on $M$ with the ${\rm weak}^*$ topology.

\bigskip
Finally let us define an exponent for the biinfinite Brownian motion.
Assume $m$ is an ergodic harmonic measure of $(M,\LL, g)$ and $h$
the characteristic harmonic function of $m$. Given $\omega\in \overline\Omega$ and a positive number $t$, 
the ratio $h(X_t(\omega))/h(X_0(\omega))$
is well defined by Proposition \ref{4}, since
a path from $X_0(\omega)$ to $X_t(\omega)$ is specified by $\omega$. 
Define a random variable $A_t:\overline\Omega\to\R$ by
$$A_t=\log h(X_t)-\log h(X_0).$$

Let us show that $A_t\in L^1(\overline \Omega,\overline\PP_m)$.
Denote the expectation of $W^x\times W_-^x$ by $\overline E^x$. Since
$A_t$ is $\FF_0$-measurable, we have $\overline E^x[A_t]=E^x[A_t]$, where $E^x$
is the expectation of $W^x$ defined before.
By the Harnack inequality 
$$
 E^x[\abs{A_t}]\leq C_1 E^x[d(X_0,X_t)] \leq C_2 t,
$$
where $d$ is the leafwise distance on the
universal cover of the leaf induced from $g$. The last inequality
follows from the bounded geometry of the leaf. 
Thus we have 
$$
\overline\EE_m\,[\abs{A_t}]=\int_ME^x[\abs{A_t}]dm(x)\leq
C_2t,
$$
showing that $A_t\in L^1(\overline \Omega, \overline \PP_m)$.

Now $A_t$ satisfies
\begin{equation} \label{e2}
A_{t+t'}=A_t+A_{t'}\circ \theta_t.
\end{equation}

This shows that $\overline\EE_m\,[A_t]$ is additive in $t$. Moreover it is
continuous in $t$ at $t=0$, since $E^x[d(X_0,X_t)]\to0$ as $t\to 0$.
That is, $\overline\EE_m[A_t]=-\lambda t$ for some number $\lambda$.

\begin{proposition} \label{10}
We have 
$\lim_{t\to\infty}(1/t)A_t=-\lambda$ almost surely, and
$\lambda\geq0$; furthermore $\lambda>0$ unless
$m$ is completely invariant.
\end{proposition}

\bd
The first statement follows from (\ref{e2}) by the Birkhoff ergodic
theorem.

To show $\lambda\geq0$ notice that
$$\int_ME^x[A_t]dm(x)=\overline\EE_m\,[A_t]=-\lambda t.$$
The expectation $E^x[A_t]$ can be computed
upstairs on the holonomy cover. Let $\hat x$ be a lift of $x$
and 
$\hat X_t(\omega)$ the lift of $X_t(\omega)$ starting at $\hat x$ 
for $\omega\in\Omega_x$.
Then
$$
E^x[A_t]=E^{x}[\log h(\hat X_t)]-\log h(\hat x)
$$
$$
\leq \log E^{x}[h(\hat X_t)]-\log h(\hat x)
= \log (\hat D_th)(\hat x)-\log h(\hat x)=0,$$
where $\hat D_t$ is the diffusion operator on the holonomy cover.
The inequality follows from the concavity of $\log$,
and the last equality from the harmonicity of $h$,
showing $\lambda\geq 0$.

For the last statement, notice that $\lambda=0$
implies that for fixed $t$, $h(\hat X_t)$ is constant 
$W^{ x}$-almost surely. This shows that $h$ is constant for the
holonomy cover of $m$-a.a.\ leaf, completing the proof.
\qed

\medskip

For $-t<0$ define a random variable $A_{-t}:\overline \Omega\to\R$ by
$$A_{-t}=\log h(X_{-t})-\log h(X_0).$$
It satisfies
\begin{equation} \label{e5}
A_{-t-t'}=A_{-t}+A_{-t'}\circ \theta_{-t}.
\end{equation}

Clearly $\overline\EE_m[A_{-t}]=\lambda$, and again by the 
Birkhoff ergodic theorem we have from (\ref{e5}):

\begin{proposition} \label{11}
$\overline\PP_m$-almost surely,
$\lim_{t\to\infty}(1/t)A_{-t}=\lambda$.
\end{proposition}

Propositions \ref{10} and \ref{11} implies
that for $m$-a.a. point $x$, we have
$W^x\times W^x_-$-almost surely
$$\lim_{t\to\infty}(1/t)A_{t}=-\lambda\ \ {\rm  and }\
\lim_{t\to\infty}(1/t)A_{-t}=\lambda,
$$
showing:

\begin{theorem} \label{12}
For a non completely invariant ergodic harmonic measure,
the characteristic harmonic function is unbounded on
the holonomy cover of $m$-a.a.\ leaf.
\end{theorem}

\section{Hyperbolic laminations}

Henceforth in this paper 
we only consider a compact hyperbolic
$C^2$ lamination $(M,\LL,g)$, i.e., 
we assume throughout that the leafwise metric $g$ has constant curvature $-1$,
and denote the dimension of leaves by $d$+$1$.
Let $m$ be an ergodic harmonic measure for $\LL$. The universal cover 
of $m$-a.a.\ leaf $L$ is identified with the simply connected complete hyperbolic
space $\DD$, and the characteristic harmonic function $h$
of $m$ is defined on $\DD$. Choose a base point $\tilde x\in\DD$
and assume $h(\tilde x)=1$.
For any point
$\xi$ of the ideal boundary $\boundary$, let $k_\xi$
denote the minimal positive harmonic function on $\DD$
corresponding to $\xi$
normalized to take value 1 at $\tilde x$.
In other words, $k_\xi=\exp(-dB_\xi)$, where 
$B_\xi$ is the Buseman function corresponding to $\xi$ such
that  $B_\xi(\tilde x)=0$. 
Then there is
a unique probability measure $\mu_{\tilde x}$ on $\boundary$ such that

\begin{equation} \label{e8}
h=\int_{\boundary} k_{\xi} d\mu_{\tilde x}(\xi).
\end{equation}

See [AS] for details and related topics.
Although the measure $\mu_{\tilde x}$ depends on the choice of
the point $\tilde x$, its equivalence class $[\mu_L]$ is an invariant
of the leaf $L$. Here two measures $\mu_1$ and $\mu_2$ on $\boundary$
are said to be equivalent if for any Borel subset $B$ of $\boundary$,
$\mu_1(B)=0$ if and only if $\mu_2(B)=0$.
In fact, for another point $\tilde y\in\DD$, we have
\begin{equation}\label{e8+}
h/h(\tilde y)=\int_{\boundary} k_{\xi}/k_{\xi}(\tilde y) d\mu_{\tilde y}(\xi).
\end{equation}
The uniqueness of the measure $\mu_{\tilde x}$ implies by (\ref{e8})
and (\ref{e8+}) that
$$\mu_{\tilde x}=(h(\tilde y)/k_\xi(\tilde y))\mu_{\tilde y},
$$
showing that $\mu_{\tilde x}$ and $\mu_{\tilde y}$ differ
by a multiple of a bounded positive function, that is,
they are equivalent.

\begin{theorem} \label{13}
For a non completely invariant
ergodic harmonic measure $m$ on a compact hyperbolic
lamination $(M,\LL,g)$ and for $m$-a.a.\ leaf $L$,
the measure class $[\mu_L]$ is singular to
the Lebesgue measure of $\boundary$.
\end{theorem}

Before starting the proof, we need to study connections
among the probability measures on $\boundary$, positive
harmonic functions on $\DD$ and the Wiener measures.

Denote by $\PPP(\boundary)$ the space of probability measures on
$\boundary$, a compact metrizable convex
space by the weak* topology. Denote by $\PPP\HH$ the space
of the positive harmonic function on $\DD$ taking value 1 at $\tilde x$,
also a compact metrizable convex space
by the compact open topology, (compact thanks to the Harnack inequality).
The map $\varphi_1:\PPP(\boundary)\to\PPP\HH$ defined by
$$
\varphi_1(\mu)=\int_{\boundary}k_\xi\,d\mu(\xi)$$
is an affine homeomorphism.

For any $f\in\PPP\HH$, define a heat kernel $q_t$ on $\DD$ by
$$
q_t(u,v)=\frac{f(v)}{f(u)}p_t(u,v),$$
where $p_t$ is the Riemannian heat kernel and $u$ and $v$
are points of $\DD$. The heat kernel defines a Wiener measure
$W^u_f$ for each point $u\in\DD$. Denote by $\Omega_{\tilde x}$
the space of continuous paths $\omega:[0,\infty)\to \DD$ such
that $\omega(0)=\tilde x$
and by $\PPP(\Omega_{\tilde x})$ the space of probability 
measures on $\Omega_{\tilde x}$. Then easy calculation shows that
the map $\varphi_2:\PPP\HH\to\PPP(\Omega_{\tilde x})$ defined by
$$
\varphi_2(f)=W^{\tilde x}_f
$$
is affine. (This is just for the base point $\tilde x$ where
$\PPP\HH$ is normalized.)

Now let $\Omega_{\tilde x}^\infty$ denotes the subspace of
$\Omega_{\tilde x}$ consisting of those paths $\omega$
in $\Omega_{\tilde x}$
such that $\lim_{t\to\infty}\omega(t)$ exists
in $\boundary$. Let us show that for any $f\in\PPP\HH$,
the set $\Omega_{\tilde x}^\infty$ is 
$W^{\tilde x}_f$-conull. As is well known, this is
true for $f=k_\xi$ for any $\xi\in\boundary$, but any
measure $W^{\tilde x}_f$ is written as the convex
integration
$$
W^{\tilde x}_f=\int_{\boundary}W^{\tilde x}_{k_\xi}d\mu(\xi)$$
for some $\mu\in\PPP(\boundary)$ since
$\varphi_1$ and $\varphi_2$ are affine, showing the claim in
the general case.

Denoting by $X_\infty:\Omega_{\tilde x}^\infty\to\boundary$
the hitting map, we define an affine map
$\varphi_3:\varphi_2(\PPP\HH)\to\PPP(\boundary)$
by $\varphi_3(W^{\tilde x}_f)=X_\infty W^{\tilde x}_f$.

Then the composite $\varphi_3\circ\varphi_2\circ\varphi_1$ is the identity 
on $\PPP(\boundary)$, since this is true for the point masses,
the map $\varphi_3\circ\varphi_2\circ\varphi_1$ is affine,
and any measure in $\PPP(\boundary)$ is a convex integral
of the point masses.

\bigskip

\noindent
{\sc Proof of Theorem \ref{13}.}\ \ Let $m$ be a non
completely invariant ergodic harmonic measure of a compact
hyperbolic lamination $(M,\LL,g)$, and let $\DD$
be the universal cover of $m$-a.a.\ leaf $L$.
A base point $\tilde x\in\DD$ is chosen and the
characteristic harmonic function $h$ normalized at 
$\tilde x$ is written as (\ref{e8}) using a probability
measure $\mu_{\tilde x}$. The Wiener measure 
$W^{\tilde x}_h$ defined by the characteristic harmonic
function $h$ corresponds to the measure $W^{\tilde x}_-$
of the reverse process in Sect.\ 3. As before denote by 
$W^{\tilde x}$ the usual Riemannian Wiener measure.
Then by Propositions \ref{10} and \ref{11}, for an appropriate
choice of $\tilde x$ we 
have $W^{\tilde x}$-almost surely 
\begin{equation}\label{ee1}
\lim_{t\to\infty}(1/t)\log(h(X_t))=-\lambda,
\end{equation}
while $W^{\tilde x}_h$-almost surely
\begin{equation}\label{ee2}
\lim_{t\to\infty}(1/t)\log(h(X_t))=\lambda,
\end{equation}
where $\lambda$ is the characteristic exponent, positive
in our case. 

On one hand, the hitting measure $X_\infty W^{\tilde x}$
of the Riemannian Wiener measure $W^{\tilde x}$ coincides with
the visible measure $\mu_0$ at $\tilde x$, which
is equivalent to the Lebesgue measure. On the other hand, the other hitting
measure $X_\infty W^{\tilde x}_h$ is the measure $\mu_{\tilde x}$.
Thus we have
$$
W^{\tilde x}=\int_{\boundary}W^{\tilde x}_{k_\xi}d\mu_0(\xi)\
\mbox{ and }\
W^{\tilde x}_h=\int_{\boundary}W^{\tilde x}_{k_\xi}d\mu_{\tilde
x}(\xi).
$$
That is, for $\mu_0$-a.a.\ point $\xi$, 
$W^{\tilde x}_{k_\xi}$-a.a.\ path satisfies 
(\ref{ee1}), while
for $\mu_{\tilde x}$-a.a.\ point $\xi$, 
$W^{\tilde x}_{k_\xi}$-a.a.\ path satisfies 
(\ref{ee2}),
showing that the two measures $\mu_0$ and $\mu_{\tilde x}$
are mutually singular.
\qed

\section{The leafwise unit tangent bundle of a hyperbolic lamination}

Associated with a compact hyperbolic lamination $(M,\LL,g)$,
there is defined
the leafwise unit tangent bundle $N$ of $\LL$ and the 
stable foliation $\HH$ on $N$. The space $M$ is covered by open sets $E_i$
on which the local charts $\varphi_i:E_i\to U_i\times Z_i$
are defined. For a hyperbolic lamination, we can assume that each
$U_i$ is an open (precompact) ball in the hyperbolic space $\DD$
and the transition function $\psi_{ji}=\varphi_j\circ\varphi_i^{-1}$
wherever defined is of the form
\begin{equation} \label{e6}
\psi_{ji}(u,z)=(g(z)u, \beta(z)),
\end{equation}
where $g(z)$ is an element of the Lie group $G$
of the orientation preserving isometries of $\DD$.
The leafwise unit tangent bunde $N$ of $\LL$
is defined from the collection of spaces $T^1(U_i)\times Z_i$
by glueing them using the transition function $\psi_{ji}$
defined by the same expression as (\ref{e6}),
where $T^1(U_i)$ is the unit tangent bundle of $U_i$. 

Notice that the tangent bundle $T^1(\DD)$ is $G$-equivariantly identified
with $\DD\times\boundary$ by assigning 
to a unit tangent vector $v$ the couple $(\pi(v),v_\infty)$,
where $\pi:T^1(\DD)\rightarrow \DD$ is the bundle projection
and $v_\infty\in\boundary$ is the hitting point of the
geodesic ray whose initial vector is $v$.

Thus a local chart $T^1(U_i)\times Z_i$ is identified with
$U_i\times\boundary \times Z_i$. Then the transition function becomes
$$
\psi_{ji}(u,\xi,z)=(g(z)u, g(z)\xi, \beta(z)).
$$
The plaques of the form $U_i\times\xi\times z$ are incorporated
to a lamination
$\HH$ of $N$, called the {\em stable foliation} of $\LL$.

The canonical projection $p:N\to M$ yields a submersion of a leaf
of $\HH$ onto a leaf of $\LL$, and thus 
the leafwise Riemannian metric $g$ of $\LL$ can be pulled up
to a leafwise Riemannian metric $\check g$ of $\HH$,
the triplet $(N,\HH,\check g)$ being a compact hyperbolic
lamination. The leafwise volume form of $\HH$ is again denoted
by ${\rm vol}$.

As before $k_\xi$ denotes the minimal positive harmonic function
associated to the point $\xi\in\boundary$ normalized at the point
$\tilde x$.

\begin{definition}
A harmonic measure $\lambda$ on $N$ is called {\em pointed
harmonic} if for each local chart $U\times\boundary\times Z$,
$\lambda$ disintegrates on a plaque $U\times\xi\times z$
to a probabality measure which is a constant times $k_\xi {\rm vol}$.
\end{definition}

The purpose of this section is to establish a one to one
correspondence between harmonic measures of $\LL$
and pointed harmonic measures of $\HH$.
We begin with a harmonic measure $m$ of $\LL$, and
associate it to a pointed harmonic measure upstairs.
Let $x$ be a point on $m$-a.a.\ leaf $L$ of $\LL$,
and let $\tilde x$ be a lift of $x$ to the universal
cover $\DD$ of $L$. Then a probability measure 
$\mu_{\tilde x}$ on $\boundary$ is defined using the characteristic
harmonic function $h$ normalized at $\tilde x$ as in (\ref{e8}).

On the other hand, the unit tangent space $T^1_xL$ is identified
with its lift $T^1_{\tilde x}\DD$, and the latter with
$\boundary$ by the visible map. By these identifications,
the measure $\mu_{\tilde x}$ on $\boundary$ corresponds
to a measure $\mu_x$ on $T^1_xL$, the notation being judtified
by the following lemma.

\begin{lemma} \label{31}
The measure $\mu_x$ is independent of the choice
of a lift $\tilde x$ of $x$.
\end{lemma}

\bd
We have only to prove that if $\gamma$ is a deck transformation
of the covering map $\DD\to L$, then 
$\mu_{\gamma \tilde x}=\gamma\mu_{\tilde x}$.
In this proof, we need a refined notation: the minimal 
positive harmonic function associated to $\xi\in\boundary$
is denoted by $k_{\xi,\tilde x}$ in order to keep in mind the
point $\tilde x$ where it is normalized.
Clearly we have
$$
k_{\gamma\xi,\gamma\tilde x}\circ\gamma=k_{\xi,\tilde x}.
$$
On the other hand by the definition of $\mu_{\tilde x}$,
the characteristic harmonic function $h$ normalized at $\tilde x$
is given by
$$
h=\int_{\boundary} k_{\xi,\tilde x}d\mu_{\tilde x}(\xi).$$
Therefore
\begin{equation} \label{e31}
h\circ\gamma^{-1}=\int_{\boundary}k_{\gamma\xi,\gamma\tilde
x}d\mu_{\tilde x}(\xi)
=\int_{\boundary}k_{\xi,\gamma\tilde x}d(\gamma\mu_{\tilde x})
(\xi).
\end{equation}
Now by Corollary \ref{4+}, 
$h\circ\gamma^{-1}$ is a constant multiple of $h$,
normalized at the point 
$\gamma \tilde x$. Therefore we have:
\begin{equation} \label{e31a}
h\circ\gamma^{-1}=\int_{\boundary}k_{\xi,\gamma\tilde
x}d\mu_{\gamma\tilde x}(\xi).
\end{equation}
Comparing (\ref{e31}) with (\ref{e31a}), 
the uniqueness of the probability measure
shows that $\mu_{\gamma\tilde x}=\gamma\mu_{\tilde x}$.
\qed

\bigskip
The inclusion $T_x^1\LL\hookrightarrow N$ induces
a map from $\PPP(T_x^1\LL)$ to $\PPP (N)$
among the spaces of the probability measures.
The image of $\mu_x$ by this map is also denoted by the same
letter, by abuse of notations.

Recall that if $(X,\mu)$ is a measured space and  $(Z,\mathcal B)$ is a
Borel space, then a map $\psi:X\to Z$ is called measurable
if for any $B\in\mathcal B$, $\psi^{-1}(B)$ is a measurable set.
Of course this depends only on the equivalence class of the measure $\mu$.
If $Z=\PPP(Y)$, the space of the probability measures of
a compact metric space $Y$, then $\psi:X\to\PPP(Y)$ is
said to be measurable if it is measurable with respect to
the Borel structure of $\PPP(Y)$ associated with the
weak* topology.
This is equivalent to saying that $x\mapsto\psi(x)(f)$
is measurable for any continuous function $f$ on $Y$.

\begin{lemma} \label{32}
The assignment $M\ni x\mapsto\mu_x\in\PPP(N)$ is measurable with respect to
the harmonic measure $m$.
\end{lemma}

\noindent
{\sc Proof.}
Since for any local chart $U\times Z$ of $\LL$, $U$
is assumed to be a domain in $\DD$, the inclusion map of $U\times Z$ into
$M$ can be extended using leafwise geodesics to a lamination preserving 
submersion $\varphi:\DD\times Z\to M$ in such a way that it
is a local isometry on each leaf. The set $\DD\times Z$ is
called a {\em prolonged local chart} of $\LL$.
Associated to it we have a prolonged local chart
$\DD\times\boundary\times Z$ for $\HH$.

By Theorem \ref{2}, the harmonic measure $m$ restricted to a local
chart $U\times Z$ is given by
$$m\vert_{U\times Z}=\int_Zh{\rm vol}\,d\nu,$$
where $h$ is a measurable function defined on $U\times Z$,
harmonic on a plaque $U\times z$ for $\nu$-a.a.\ $z$.
For the prolonged
local chart $\DD\times Z$ let $m\vert_{\DD\times Z}$
be the {\em lift} of $m$ to $\DD\times Z$, i.e., the integral
of the counting measure on the fiber of the submersion
$\DD\times Z\to M$ over $m$.
Then we still have
\begin{equation} \label{e32}
m\vert_{\DD\times Z}=\int_Z h\, {\rm vol} d\nu,
\end{equation}
where $h$ is an obvious extension.
Notice that a slight generalization of Theorem \ref{2} shows that 
$h$ is measurable with respect to $m\vert_{\DD\times Z}$.

Denote by $\PPP\HH_u$ the space of positive harmonic functions
taking value 1 at $u\in\DD$. Then there is an affine homeomorphism
of $\PPP\HH_u$ with $\PPP(\boundary)$. 
Let $(u,z)\in \DD\times Z$ corresponds to $x\in M$ by the submersion.
The measure $\mu_x=\mu_{(u, z)}$ of $\boundary$ is associated
to the function $h(\cdot,z)/h(u,z)\in\PPP\HH_u$ by the above homeomorphism.

\begin{sublemma} \label{32a}
The assignment 
$\DD\times Z\ni(u,z)\mapsto\mu_{(u,z)}\in\PPP(\boundary)$
is measurable with respect to $m\vert_{\DD\times Z}$.
\end{sublemma}

\bd
The measure $m\vert_{\DD\times Z}$ is equivalent to ${\rm vol}\otimes\nu$.
Therefore by Fubini,
there is a vol-conull subset $\DD_*$  
such that for any poin $u\in\DD_*$, the set 
$\{z\in Z;\,h(u,z)<\alpha\}$
is $\nu$-measurable for any $\alpha\in\mathbb Q$.
It is routine to show then for any any $u\in\DD_*$
and $\alpha\in\R$,
the set $\{z\in Z;\,h(u,z)<\alpha\}$
is $\nu$-measurable. 

For any $u\in\DD_*$, the assignment to $z\in Z$ of the
harmonic function
$h(\cdot,z)/h(u,z)$ in $ \mathcal P\HH_u$
is $\nu$-measurable with respect to the $\sigma$-algebra 
$\mathcal B(\PPP\HH_u)$ of the
pointwise convergence topology on $\DD_*$. 
In fact for any $v\in\DD_*$ and $a>0$, the set
$$
\{z\in Z;\, h(v,z)>ah(u,z)\}=
\bigcup_{\alpha\in\mathbb Q}(\{z\in Z;\,h(v,z)\geq\alpha\}\cap\{
z\in Z; ah(u,z)<\alpha\})$$
is $\nu$-measurable. 

The $\sigma$-algebra $\mathcal B(\PPP\HH_u)$ 
coincides with the $\sigma$-algebra of the compact open topology.
In fact for $(a,b)\subset\R$ and a compact ball $D$ of $\DD$,
the set 
$$ \PPP\HH_u(D,(a,b))=\{f\in\PPP\HH_u;\, f(D)\subset(a,b)\}
$$
belongs to  $\mathcal B(\PPP\HH_u)$, since
for a subset $\{u_j\}_{j\in\N}\subset \DD_*\cap D$ dense in $D$, we have
\begin{eqnarray*}
&\PPP\HH_u(D,(a,b))=\bigcup_{n\in\N}\{f\in\PPP\HH_u;\, f(D)
\subset[a+n^{-1},b-n^{-1}]\}\\
&=\bigcup_{n\in\N}\bigcap_j\{f\in\PPP\HH_u;\, f(u_j)
\in [a+n^{-1},b-n^{-1}]\},
\end{eqnarray*}
and this subset belongs to $\mathcal B(\PPP\HH_u)$.
A general compact subset $K\subset\DD$ can be written
as the decreasing intersection of finite unions of compact balls $D_n$,
and the like defined  set $\PPP\HH_u(K,(a,b))$ also belongs
to $\mathcal B(\PPP\HH_u)$, since
$$\PPP\HH_u(K,(a,b))=\bigcup_n\PPP\HH_u(D_n,(a,b)).$$

The space $\mathcal P\HH_u$ with the compact open topology is homeomorphic
to the space $\PPP(\boundary)$ with the ${\rm weak}^*$ topology.
This shows the $\nu$-measurability
of $\mu_{(u,z)}$  
in the variable $z$ for any $u\in\DD_*$. On the other hand
the measure $\mu_{(u,z)}$ is continuous in the variable $u$
for any $z\in Z$.

Let $f:\boundary\to\R$ be an arbitrary continuous function
and fix it for a while.
For any $a\in\R$, define
$$
S(a)=\{(u,z)\in U\times Z;\, \mu_{(u,z)}(f)\geq a\}.$$
The proof of the sublemma 
is complete if we show that $S(a)$ is a measurable set.

For any $z\in Z$, define the $z$-slice $S(a)_z\subset U$
by 
$$
S(a)\cap(U\times z)=S(a)_z\times z.$$ 
Similarly define the $u$-slice $S(a)_{u}\subset Z$ 
for any $u\in U$.
Then $S(a)_{u}$ is $\nu$-measurable for any $u\in\DD_*$, and $S(a)_{z}$ is
closed for any $z\in Z$. 
Moreover since $\mu_{(u,z)}(f)$ is a continuous function
of $u$,   
$\{S(a_k)_z\}$ forms a (closed) neighbourhood system of $\{S(a)_z\}$
for a sequence $a_k\uparrow a$.
Choose a compact ball $D\subset \DD$ and define
$$S(a)_D=\{z\in Z;\, S(a)_z\cap D\neq\emptyset\}.
$$
Then $S(a)_D$ is $\nu$-measurable.
In fact, for a dense subset $\{u_j\}_{j\in\N}$ of $D\cap\DD_*$,
we have
$$
S(a)_D=\{z\in Z;\, S(a_k)_z\cap\{u_j\}\neq\emptyset,\ \forall k\}
=\bigcap_{k}\bigcup_j S(a_k)_{u_j}.
$$

Now let $\{\mathcal D_i\}$ be a sequence of coverings of $\DD$ by 
countably many compact
balls such that ${\rm mesh}(\mathcal D_i)\to0$ as $i\to\infty$. 
Define
$$
S(a)^i=\bigcup_{D\in\mathcal D_i}D\times S(a)_D.$$
Then $S(a)^i$ is measurable.
On the other hand, since $S(a)_z$ is closed, we have
$S(a)_z=\bigcap_iS(a)_z^i$. That is, $S(a)=\bigcap_iS(a)^i$
and $S(a)$ is measurable,
completing the proof.
\qed

\bigskip

\noindent
Sublemma \ref{32a} implies in particular for any local chart $U\times Z$,
the assignment 
$$
U\times Z\ni (u,z)\mapsto \mu_{(u,z)}\in\PPP(\boundary)
$$
is measurable.  
Define a map $\iota_{(u,z)}:\boundary\to U\times\boundary\times Z$ 
by $\iota_{(u,z)}(\xi)=(u,\xi,z)$.
Consider a map
$$
\psi:U\times Z\times \PPP(\boundary)\to\PPP(U\times \boundary\times Z)
$$
defined by $\psi(u,z,\mu)=\iota_{(u,z)}\mu$. 
Consider also a map $\phi:\PPP(U\times \boundary\times Z)\to\PPP(N)$ induced
by the inclusion.
If $(u,z)\in U\times Z$ corresponds to $x\in M$, then
\begin{equation} \label{mu}
\mu_x=(\phi\circ\psi)(u,z,\mu_{(u,z)}).
\end{equation}

The proof of Lemma \ref{32} is complete if we show that the RHS
of (\ref{mu}) is a measurable function of $(u,z)$.
Here we have:
\begin{sublemma} \label{32c}
The map
$\psi:U\times Z\times\PPP(\boundary)\to\PPP(U\times \boundary\times Z)$
is continuous.
\end{sublemma}

\bd
Denote by $C_c(U\times\boundary\times Z)$ the
space of the continuous functions with compact supports.
Consider a product $f\circ p_1\,g\circ p_2$
($f\in C_c(U\times Z)$, $g\in C(\boundary)$,
$p_1:U\times\boundary\times Z\to U\times Z$, 
$p_2:U\times\boundary\times Z\to \boundary$, the canonical projections).
Then clearly
$$
\psi(u,v,\mu)(f\circ p_1 g\circ p_2)=f(u,v)\mu(g)$$
is a continuous function of $(u,v,\mu)$.
On the other hand, finite sums of the products $f\circ p_1\,g\circ p_2$
form a dense subset of $C_c(U\times \boundary\times Z)$
in the topology of the uniform convergence on compact sets.
Standard argument shows that $\psi(u,v,\mu)(F)$ is continuous for
any $F\in C_c(U\times\boundary\times Z)$, finishing the proof.
\qed

\bigskip
On the other hand, $\phi$ is obviously continuous.
The RHS of (\ref{mu}) is now shown to be a measurable function
of $(u,z)$,
completing the proof of Lemma \ref{32}.
\qed

\bigskip

Integrating the measurable system of probability measures 
$\{\mu_x\}_{x\in M}$
over $m$, we obtain a probability measure $\lambda(m)$ of $N$, called
the {\em canonical lift} of $m$.

\begin{theorem} \label{33}
For any harmonic measure $m$ of $M$,
the canonical lift $\lambda(m)$ is pointed harmonic.
\end{theorem}

\bd
Recall that on a prolonged local chart $\DD\times Z$, the lift of the
harmonic measure $m$ is written as in (\ref{e32}), and 
the canonical lift $\lambda(m)$ on the associated prolonged
local chart $\DD\times \boundary\times Z$ disintegrates
on $\DD\times\boundary\times z$ to (a constant multiple of)
\begin{equation} \label{e34}
\int_{\DD}h(u)\mu_u\,d{\rm vol}(u),
\end{equation}
where $\mu_u$ is the probability measure in $\PPP(\boundary)$
determined by the equality
$$
\frac{h(v)}{h(u)}=\int_{\boundary}\frac{k_\xi(v)}{k_\xi(u)}d\mu_u(\xi),
\ \ \forall v\in\DD,$$
where
$k_\xi$ is the minimal harmonic function normalized
at the base point $\tilde x$.

In order to disintegrate further the measure in (\ref{e34}) 
on $\DD\times\xi\times z$,
we have to transform the measure $\mu_u$ which
depends on $u\in\DD$ to a fixed measure $\mu_{\tilde x}$.
First of all we have
$$
h(v)=\int_{\boundary}h(u)\frac{k_\xi(v)}{k_\xi(u)}d\mu_u(\xi)
=\int_{\boundary}k_\xi(v)\frac{h(u)}{k_\xi(u)}\frac{d\mu_u}
{d\mu_{\tilde x}}(\xi)d\mu_{\tilde x}(\xi).
$$

Hence by the uniqueness of the probability measure, we have
$$
\frac{h(u)}{k_\xi(u)}\frac{d\mu_u}{d\mu_{\tilde x}}(\xi)=1,$$
showing that
$$
\int_{\DD}h(u)\mu_u\,d{\rm vol}(u)
=\int_{\DD}k_\xi(u)\mu_{\tilde x}\,d{\rm vol}(u).$$
This implies that the lift of the measure $\lambda(m)$
disintegrates on $\DD\times \xi\times z$ to a constant
multiple of $k_\xi\,{\rm vol}$, completing the proof.
\qed

\bigskip

Conversely given any pointed harmonic measure on the leafwise
unit tangent space $N$, 
its push down is a harmonic measure on
$M$ by Theorem \ref{2}. It is easy to show the following
theorem by analogous computation.

\begin{theorem} \label{15}
A harmonic measure on a compact hyperbolic
lamination $(M,\LL,g)$ corresponds one to one to a
pointed harmonic measure on its leafwise unit
tangent bundle $(N,\HH,\check g)$, by the operations of
taking the canonical lift
and pushing down.
\end{theorem}

\begin{example} \label{16}
If $M$ is a closed oriented hyperbolic surface,
considered as a single leaf lamination, then the unique
harmonic measure $m$ is the (normalized) area form. 
The canonical lift $\lambda(m)$ on the
unit tangent bundle $T^1M$ is the (normalized) Haar measure.
\end{example}

\begin{remark}
In case $d=1$ the minimal parabolic subgroup $P$ of $G$
acts on the leafwise tangent bundle $N$ 
of a compact 2 dimensional hyperbolic lamination from the right
in such a way that the orbit lamination
is the stable foliation $\HH$,
and a probability measure of $N$ is pointed harmonic if and only
if it is invariant by the action of $P$.
Theorem \ref{15} in this case is already obtained in [M]
and [BM] by a somewhat different dynamical method.
For higher dimension we do not have such description of
pointed harmonic measures.
\end{remark}

\section{The dichotomy}

Let $m$ be a harmonic measure on a compact hyperbolic lamination
$(M,\LL,g)$. For $m$-a.a.\ leaf $L$, we have defined a measure
class $[\mu_L]$ on the boundary $\boundary$ of the universal cover
$\DD$ of $L$. In this section we shall prove Theorem \ref{1},
i.e.,  that for an ergodic harmonic measure $m$, either the support
$K_L={\rm Supp}([\mu_L])$, called the {\em characteristic set}
of $L$, is a singleton for any $m$-a.a.\ leaf,
or is the total space $\boundary$. 

The argument closely follows the proof of Proposition 1 [MV],
in which the authors attribute the original idea to Etienne Ghys.

To begin with, let us notice the following fact.
Let $\Gamma$ be the group of deck transformations of the
covering map $\DD\to L$. In the proof of Lemma \ref{31}, 
we have shown that $\mu_{\gamma \tilde x}=\gamma\mu_{\tilde x}$
for any $\gamma\in\Gamma$.
On the other hand the equivalence class
of the measure $\mu_{\tilde x}$ does not depend on the choice of
the particular point $\tilde x$ from $\DD$, as is explained in the
beginning of Sect.\ 4.
This shows that $\gamma K_L=K_L$. 

Given a closed subset $K$ of $\boundary$ which is not a singleton,
the {\em convex hull} of $K$, denoted by $C(K)$,
is the convex hull in $\DD$ of the union of 
all the geodesics joining two points
of $K$.
It is a closed convex subset of $\DD$, and the assignment
$K\mapsto C(K)$ is $G$-equivariant, where $G$ is the group
of all the orientation preserving isometries of $\DD$.
Therefore we have:

\begin{lemma} \label{17}
Assume $K_L$ is not a singleton. Then the convex hull $C(K_L)$
of $K_L$, as well as its closed $r$-neighbourhood $N_L(r)$ {\rm ($r>0$)}, 
is a $\Gamma$-invariant subset of $\DD$.
\end{lemma}

Choose a prolonged local chart $\DD\times Z$, and denote
the characteristic set of the leaf of $\LL$ corresponding
to $\DD\times z$ by $K_z$. Denote by $\CC(\boundary)$
the space of nonempty closed subsets of $\boundary$,
equipped with the $\sigma$-algebra ${\mathcal B}_{\mathcal C}$ of
the Hausdorff topology.

\begin{lemma} \label{18}
The assignment $Z\ni z\mapsto K_z\in\CC(\boundary)$ is $\nu$-measurable with respect to
${\mathcal B}_{\mathcal C}$.
\end{lemma}

\bd
For any open
subset $U$ of $\boundary$, define ${\mathcal C}(\boundary)_U$
to be the open subset of $\CC(\boundary)$ consisting of those
closed sets which intersects $U$. It is well known, easy to show,
that the open sets $\CC(\boundary)_U$ 
generate the $\sigma$-algebra ${\mathcal B}_{\mathcal C}$.
Therefore it suffices to show that the set 
$$
Z_U=\{z\in Z;\, K_z\in \CC(\boundary)_U\}$$
is $\nu$-measurable. Choose a countable family $\{f_i\}$
of nonnegtive continuous functions supported in $U$ such that
the union of their support is $U$, and take a base point
$\tilde x\in\DD_*$, where $\DD_*$ is the subset defined in 
the proof of Sublemma \ref{32a}. Then 
the set $Z_U$ consists of exactly those points $z$
such that $\mu_{(\tilde x,z)}(f_i)>0$ for some $i$.
The $\nu$-measurable dependence of $\mu_{(\tilde x,z)}$ 
in the variable $z$
established in the proof of Sublemma \ref{32a} completes the proof.
\qed

\begin{definition}
(1) Let $M_{\rm I}$ be the union of $m$-a.a.\ 
leaves $L$ such that the characteristic set $K_L$ is a singleton.

(2) Let $M_{\rm II}$ be the union of $m$-a.a.\ leaves $L$ such that
$K_L=\boundary$.

(3) Let $M_{\rm III}=M\setminus(M_{\rm I}\cup M_{\rm II})$.
\end{definition}

Lemma \ref{18} implies that the three subsets are $m$-measurable.
Since they are saturated and the harmonic measure $m$ is
ergodic, only one of them is conull. 
Henceforth we assume that $M_{\rm III}$ is conull and deduce a
contradiction,
which is sufficient for the proof of Theorem \ref{1}.
For any $m$-a.a.\ leaf $L$ and for $r>0$, 
consider the image of $N_L(r)$ by the covering map $\DD\to L$.
Taking their union for any $m$-a.a.\ leaf $L$,
we get a subset of $M$, denoted by $N(r)$.

\begin{lemma} \label{19}
The subset $N(r)$ is measurable.
\end{lemma}

\bd
Denote by $\CC(\DD\cup\boundary)$ the set of nonempty
closed subsets of the compactification $\DD\cup\boundary$,
equipped with the Hausdorff topology.
Then the map from $\CC(\boundary)$ to $\CC(\DD\cup\boundary)$ 
which assigns to $K$ the closure 
of the $r$-neighbourhood of the convex hull
of $K$ is clearly continuous.

Consider a prolonged local chart $\DD\times Z$ and again let $K_z$
denote the characteristic set of the leaf in $\LL$ which
corresponds to $\DD\times z$. Also denote by $N_z(r)\subset\DD$
the closed $r$-neighbourhood of the convex hull of $K_z$.
Then by the above observation and by Lemma \ref{18},
the map
$$Z\ni z\mapsto N_z(r)\cup K_z\in\CC(\DD\cup\boundary)$$
is measurable. In particular for any
open subset $U$ of $\DD$, the set
$$
\{z\in Z;\, N_z(r)\cap U\neq\emptyset\}
$$
 is a measurable
subset of $Z$.

Let us show that the union $N_Z(r)=\bigcup_z N_z(r)\times z$ is a measurable
subset of $\DD\times Z$. Choose a sequence of open coverings
of $\DD$,
${\mathcal U}_1\prec{\mathcal U}_2\prec\cdots$, such that
mesh$({\mathcal U}_i)\to 0$. Define
$$
N_Z(r)^i=\bigcup_{U\in{\mathcal U}_i}U\times\{z;\, N_z(r)
\cap U \neq\emptyset\}.
$$
Then the set $N_Z(r)^i$ is measurable, and hence
$N_Z(r)=\bigcap_iN_Z(r)^i$ is also measurable.

Now the image of $N_Z(r)$ by the submersion of $\DD\times Z$ to $M$ 
is measurable.
In fact, $N_Z(r)$ is a union of a null set and a 
Borel set. The image of a null set is null by the
definition of the lift $m\vert_{\DD\times Z}$ of $m$. 
On the other hand the image
of a Borel set by a countable to one Borel map is Borel.
This is a well known fact about standard Borel spaces,
and follows e.g. from [Ke] Corollary 15.2 and [S]
Theorem 1.3.
Now the set $N(r)\subset M$ is a finite union of measurable sets
and is measurable.
\qed

\bigskip

Let us finish the proof of Theorem \ref{1}. 
Recall we are assuming that $M_{\rm III}$ is conull in way of
contradiction.
Since
$M=\bigcup_rN(r)$ mod 0, we have $m(N(r))>0$ for some $r$.
By Theorem \ref{5}, the measure $\PP_m$ on the space $\Omega$
of leafwise paths is ergodic with respect to
the shift semiflow $\theta_t$. This means that for $\PP_m$-almost any
path the average time of stay in the set $X_0^{-1}(N(r))$ 
is equal to $\PP_m(X_0^{-1}(N(r)))=m(N(r))$.
In other words for $m$-a.a.\ $x$, $W^x$-almost
surely we have
\begin{equation} \label{e10}
\lim_{t\to\infty}\frac{1}{t}dt\{s\in[0,t] ;\, X_s\in N(r)\}=m(N(r))>0,
\end{equation}
where $dt$ denotes the Lebesgue measure on $[0,\infty)$.

But by Lemma \ref{17}, the inverse image
$p^{-1}(N(r))$ of the universal covering map $p:\DD\to L$ 
of $m$-a.a.\ leaf $L$ coincides with
the set $N_L(r)$, the closed $r$-neighbourhood of the
convex hull of the characteristic set $K_L$.
Since $K_L\neq\boundary$, there is a closed
nondegenerate interval $I$ contained in $\boundary\setminus K_L$.
For any point $x$ on $L$, the set of
paths whose lifts hit $I$ has positive $W^x$-measure.
On the other hand for those paths the limit of (\ref{e10}) must be 0,
since there is a neighbourhood of $I$ 
in $\DD\cup\boundary$ which does not intersect $N(r)$. A contraction.
Theorem \ref{1} is now proved.

\bigskip

\begin{example}
For any harmonic measure $m$ of a compact hyperbolic lamination,
the canonical lift $\lambda(m)$ of $m$, a pointed
harmonic measure of the leafwise
unit tangent bundle, is of type I. Especially
the unique ([G], [DK]) harmonic measure of
the Anosov foliation on
the unit tangent bundle of a closed oriented hyperbolic surface
is of type I. 
\end{example}

Ergodic completely invariant measures are typical examples of 
harmonic measures of Type II. But there are some more.
An example is in order. Let
$\Sigma=\Gamma\setminus{\mathbb D}^2$, where 
$\Gamma<PSL(2,\R)$ is a purely hyperbolic cocompact
Fuchsian group. 

Choose any homomorphism $\rho:\Gamma\to {\rm Homeo}(Z)$
to the group of the homeomorphisms of a compact metric space $Z$
which satisfies the following conditions.

(1) The homomorphism $\rho$ is not injective.

(2) There is no $\rho(\Gamma)$-invariant measure on $Z$.

Let $M=\Gamma\setminus({\mathbb D}^2\times Z$), where the action of $\Gamma$
is by deck transformation on the first factor and by $\rho$ on
the second. Then the horizontal lamination 
$\{{\mathbb D}^2\times z\}$ on ${\mathbb D}^2\times Z$
induces a lamination $\LL$ on $M$, called the {\em suspension}
of $\rho$.
Let $m$ be any ergodic harmonic measure of $\LL$, and notice that
$m$ is 
not completely invariant by (2). 

\begin{proposition}
The above ergodic harmonic measure $m$ is of type II.
\end{proposition}

\bd
By Theorem \ref{2}, the harmonic measure $m$
determines the class of a probability measure $\nu$
on $Z$. The measure $\nu$ is quasi-invariant by the action
of $\rho(\Gamma)$.

Assume for contradiction that $m$ is of type I. Then
for the prolonged local chart ${\mathbb D}^2\times Z$,
the charcteristic set $K_z$ ($z\in Z$) is a singleton
for $\nu$-a.a.\ $z$. This yields a measurable map
$k: Z\to\oboundary$, by Lemma \ref{18}. The map $k$
is $\Gamma$-equivariant with respect to $\rho$ and the Fuchsian group
action on $\oboundary$, i.e.,  we have
$$
k(\rho(\gamma)z)=\gamma k(z)\ \ \mbox{ for all }\ \gamma\in\Gamma,
\ \ \nu{\rm -a.a.}\ \ z\in Z.
$$

The push forward measure
$k\nu$ is kept quasi-invariant by the Fuchsian group, and 
in particular its support is an infinite set.
Choose  a nontrivial $\gamma\in\Gamma$ from the kernel of $\rho$,
and let $F$ be a Borel fundamental domain of $\gamma$ for
its action on
$\oboundary\setminus{\rm Fix}(\gamma)$.  Then we have $\nu(k^{-1}(F))>0$.
On the other hand we have
$$
k^{-1}\gamma F=\rho(\gamma^{-1})k^{-1}\gamma F=k^{-1}F\ \ {\rm mod}\ 0.
$$
Thus we have $\nu(\emptyset)=\nu(k^{-1}F\cap k^{-1}\gamma F)>0$.
A contradiction.
\qed

Finally let us pose some problems.

\begin{question} It is known [K2] that 
a compact hyperbolic lamination with a type I ergodic 
harmonic measure is an amenable
measured foliation in the sense of [AR]. 
Is the converse true?
\end{question}

\begin{question}
For an ergodic harmonic measure of type I of a compact
hyperbolic lamination of dimension $d$+$1$, the characteristic
exponent satisfies $\lambda=d^2$. Is it true for type
II measure that $\lambda<d^2$?
\end{question}

\begin{question}
For an injective homomorphism from $\Gamma$ (as above) 
to $PSL(2,\R)$ with dense image, 
is the harmonic measure
of the suspension foliation  type II?
\end{question}

\bigskip
\noindent
Department of Mathematics
\\
College of Science and Technology
\\
Nihon University
\\
1-8-14 Kanda, Surugadai
\\
Chiyoda-ku, Tokyo, 101-8308 
\\
Japan

\medskip
\noindent
{\em E-mail address}: matsumo@math.cst.nihon-u.ac.jp

\end{document}